\documentclass{amsart}[11pt]

\usepackage{amscd,amssymb,latexsym,url,verbatim,graphicx,color}
\usepackage{hyperref,cleveref}
\usepackage{tikz,tikz-cd}

\usepackage[symbol]{footmisc}

\usetikzlibrary{decorations.pathmorphing}

\usepackage{cases,amsmath}

\usepackage{pxfonts}
\usepackage[euler-digits,small]{eulervm}

\usepackage{tikz,tikz-cd}

\usepackage[dvips]{epsfig}

\title{Stirling complexes}

\author{Dmitry N. Kozlov \\
\MakeLowercase{with an appendix by} \smaller{Roy Meshulam}}
\address{Department of Mathematics, University of Bremen, 28334
  Bremen, Federal Republic of Germany.}
	
\email{dfk@math.uni-bremen.de}

\address{Okinawa Institute of Science and Technology Graduate University,
1919-1 Tancha, Onna-son, Kunigami-gun,
Okinawa, Japan.}


\newtheorem{theorem}{Theorem}[section]
\newtheorem{df}[theorem]{Definition}
\newtheorem{thm}[theorem]{Theorem}

\newtheorem{prop}[theorem]{Proposition}

\newtheorem{openq}{Open Question}
\newtheorem{rem}[theorem]{Remark}

\newcommand{\nin}{\noindent}
\newcommand{\pr}{\nin{\bf Proof.} }
\newcommand{\prn}[1]{\nin{\bf Proof of #1.} }

\newcommand{\cld}{{\mathcal D}}

\newcommand{\cle}{{\mathcal E}}

\newcommand{\cls}{{\mathcal S}}

\newcommand{\rr}{{\mathbb R}}
\newcommand{\zz}{{\mathbb Z}}

\newcommand{\ra}{\rightarrow}
\newcommand{\sm}{\setminus}
\newcommand{\supp}{\text{\rm supp}\,}

\newcommand{\id}{\textrm{id}}

\newcommand{\sph}{{\mathbb S}}
\newcommand{\str}{\mathcal{S}tr}

\newcommand{\mycap}[1]{\caption{#1}}

\newcommand{\myfig}[2]{
\begin{figure}[hbt]
\centering
\input{#1.pdf_t}
\mycap{#2.}
\label{fig:#1}
\end{figure}}

\newcommand{\myset}[2]{\{#1\,|\,#2\}}

\newcommand{\demph}{\bf}

\newcommand{\ctop}{\mathbf{Top}}
\newcommand{\hocolim}{{\tt hocolim}\,}
\newcommand{\susp}{{\tt susp}\,}

\graphicspath{{./figs/}}

\numberwithin{equation}{section}
\numberwithin{figure}{section}
\numberwithin{table}{section}

\newcommand{\sn}[2]{\genfrac{\{}{\}}{0pt}1{#2}{#1}}
\newcommand{\snd}[2]{\genfrac{\{}{\}}{0pt}0{#2}{#1}}
\newcommand{\stc}[2]{\mathcal{S}tr(#1,#2)}
\newcommand{\stm}[3]{\mathcal{S}tr(#1,#2,#3)}
\newcommand{\srf}[2]{{\tt{SF}}(#1,#2)}
\newcommand{\mynewpage}{}

\begin{document}

\begin{abstract}
In this paper we study natural reconfiguration spaces associated to
the problem of distributing a fixed number of resources to labeled
nodes of a tree network, so that no node is left empty.  These spaces
turn out to be cubical complexes, which can be thought of as
higher-dimensional geometric extensions of the combinatorial Stirling
problem of partitioning a set of named objects into non-empty labeled
parts.

As our main result, we prove that these Stirling complexes are always
homotopy equivalent to wedges of spheres of the same
dimension. Furthermore, we provide several combinatorial formulae to
count these spheres.

Somewhat surprisingly, the homotopy type of the Stirling complexes
turns out to depend only on the number of resources and the number of
the labeled nodes, not on the actual structure of the tree network.
\end{abstract}

\maketitle

\section{Stirling complexes}

\subsection{Motivation}
Consider the situation where $n$ unique resources need to be
distributed among $m$ locations. Clearly, subject to the only
condition that $n\geq m$, this can be done in many different
ways. Specifically, the number of solutions is equal to $m!\sn{m}{n}$,
where $\sn{m}{n}$ is the \emph{Stirling number of the second kind},
which is a classical combinatorial function, counting the number of
ways $n$ objects can be partitioned into $m$ non-empty groups, see
\cite{GKP,Kn,S}.

Imagine furthermore, that the locations, to which the resources are
distributed, are connected by a tree network, and that each resource
can be shifted from its location to a neighboring one. Simultaneous
multiple shifts of different resources are allowed, as long as at any
point of this shifting procedure there remain some resources, which
are not being moved, in each node. We would like to model this
situation by introducing a higher-dimensional parameter space which
encodes the interplay of such shifts. In what follows we introduce a
family of combinatorial cubical complexes, which fulfill this task. We
shall call these complexes the \emph{Stirling complexes}.

In recent years topology has increasingly been used in applications, most
notably in data analysis, see \cite{Ca} and the references therein.
The idea of using higher-dimensional cell complexes to record
transformations of combinatorial objects has been a further major thread in
the tapestry of applied topology.  For instance, a family of
prodsimplicial complexes has been constructed in \cite{BaK06}, see
also \cite{BaK07,Ko07,Ko08}, to find topological obstructions to graph
colorings, a famously notorious problem. Another example is provided
by the so-called protocol complexes, which have been introduced as a
part of the topological approach to questions in theoretical
distributed computing, see \cite{HKR} and the numerous references
therein.  Optimally, such constructions provide deeper insight into
the original combinatorial questions, yielding at the same time
interesting, often highly symmetric families of combinatorial cell
complexes.

In what follows, we shall use standard facts and terminology of graph
theory, as well as algebraic topology. If the need arises, the reader
is invited to consult \cite{Har} for graph theory, and
\cite{FFG,Fu,GH,Hat,Ko08,Ko20,Mu} for algebraic topology.

%\mynewpage

\subsection{Definition of the Stirling complexes}
Let $m$ be an arbitrary integer, $m\geq 2$, and let $T$ be an
arbitrary tree on $m$ vertices, labeled with numbers $1$ through
$m$. This tree models our network.

Assume furthermore we have $n\geq m$.  We can view $T$ as a
$1$-dimensional simplicial complex, which leads us to considering the
cubical complex $T^n$.  Let us make the following observations about
this complex.
\begin{itemize}
\item The cubes of $T^n$ are indexed by the $n$-tuples
  $c=(c_1,\dots,c_n)$, where each $c_i$ is either a vertex or an edge
  of $T$.
\item The dimension of $c$ is equal to the number of $c_i$'s which are
  edges. Accordingly, the vertices of $T^n$ are indexed by the
  $n$-tuples of the vertices of $T$, the dimension of $T^n$ is equal
  to $n$, and the top-dimensional cubes are indexed by the $n$-tuples
  of the edges.
\item The boundary cubes of $c$ are obtained by replacing edges in the
  indexing $n$-tuple with adjacent vertices. The number of replaced
  edges is precisely the codimension of the corresponding boundary
  cube.
\end{itemize}

We are now ready to define our main objects of study.

\begin{df} \label{df:stc}
Given a tree $T$ with $m\geq 2$ vertices, and a positive integer $n$,
the {\demph Stirling complex} $\stc{T}{n}$ is the subcomplex of $T^n$
consisting of all $n$-tuples $c=(c_1,\dots,c_n)$, such that each
vertex of $T$ occurs as an entry in $c$ at least once.
\end{df}

Since the condition of \Cref{df:stc} is preserved by taking the
boundary, the Stirling complexes are well-defined. The following facts
hold for Stirling complexes.

\begin{itemize}
\item If $n<m$, the condition in \Cref{df:stc} cannot be fulfilled,
so $\stc{T}{n}$ is empty in this case.
\item The complex $\stc{T}{m}$ consists of $m!$ vertices, indexed by
  all permutations of the set $[m]=\{1,\dots,m\}$.
\item In general, the vertices of $\stc{T}{n}$ are indexed by all ways
  to partition the set $[n]$ into $m$ labeled parts. Accordingly, the
  number of vertices of $\stc{T}{n}$ is equal to $m!\sn{m}{n}$.
\item The dimension of $\stc{T}{n}$ is equal to $n-m$, since this is the maximal number
of resources which can be assigned to the edges of~$T$.
\end{itemize}

For each $0\leq d\leq n-m$, the Stirling complex $\stc Tn$ has
$\binom{n}{d}(m-1)^dm!\sn{m}{n-d}$ cubes of dimension $d$.  To see
this, first choose $d$ resources among $n$, then assign each resource
to one of the $m-1$ edges, and then finally distribute the rest of the
resources to the nodes, so that no node is left empty.  This gives us
the following formula for the Euler characteristic:
\begin{equation} \label{eq:chistc}
\chi(\stc Tn)=\sum_{d=0}^{n-m}(-1)^d\binom{n}{d}(m-1)^dm!\snd{m}{n-d}.
\end{equation}
\nin In what follows, we shall derive a better formula for $\chi(\stc Tn)$.

\mynewpage

\subsection{Examples} To acquaint ourselves with the Stirling complexes, let us consider
a few further examples.

\vspace{1ex}

\nin
{\it Example 1.}
The first interesting example is $\stc{T}{m+1}$.  The dimension of
this Stirling complex is 1, so it is a graph.  The numerical data of
this graph is the following.
\begin{itemize}
\item The number of vertices of $\stc{T}{m+1}$ is
\[m!\snd{m}{m+1}=m!\binom{m+1}{2}=\frac{m}{2}(m+1)!.\]
The vertices  of $\stc{T}{m+1}$ are indexed by the $(m+1)$-tuples of the vertices of $T$,
with one vertex repeating twice and all other vertices occurring
exactly once.
\item As a graph $\stc{T}{m+1}$ has $(m-1)(m+1)!$ edges; the edges are indexed by
  $(m+1)$-tuples consisting of one edge and $m$ vertices of $T$, with
  each vertex repeating exactly once.
\end{itemize}
Accordingly, the Euler characteristic of this Stirling complex is given by
\[\chi(\stc{T}{m+1})=-\dfrac{1}{2}(m-2)(m+1)!.\]

It is easy to see using a direct argument that $\stc{T}{m+1}$ is
always connected.  Therefore it is homotopy equivalent to a wedge of
$\dfrac{1}{2}(m-2)(m+1)!+1$ circles.

Consider now the special case when $m=4$. Let $T_1$ be the tree with
one vertex of degree $3$ and $3$ leaves. Let $T_2$ be the string with
$3$ edges: it has $2$ vertices of degree $2$ and $2$ leaves.  Both
$\stc{T_1}{5}$ and $\stc{T_2}{5}$ are connected and have $240$
vertices and $360$ edges.  However, these two graphs are different:
$\stc{T_1}{5}$ has $60$ vertices with valency $6$, and the rest of the
vertices with valency $2$, whereas all vertices of $\stc{T_2}{5}$ have
valency $2$ or $4$. We see therefore that, while topology of
$\stc{T_1}{5}$ and $\stc{T_2}{5}$ is the same, the spaces themselves
depend on the actual tree structure of $T_1$ and~$T_2$.

%\mynewpage

\vspace{1ex}

\nin
{\it Example 2.}
Next, consider the cubical complexes $\stc{T}{m+2}$, for $m\geq 2$. These
are $2$-dimensional.
The number of vertices is given by
\[f_0:=m!\sn{m}{m+2}=m!\frac{1}{24}m(m+1)(m+2)(3m+1)=
m(3m+1)\frac{(m+2)!}{24}.\]

The number of edges is given by
\[f_1:=(m+2)(m-1)\frac{m}{2}(m+1)!=12m(m-1)\frac{(m+2)!}{24}.\]

Finally, the number of squares is given by
\[f_2:=\binom{m+2}{2}(m-1)^2m!=12(m-1)^2\frac{(m+2)!}{24}.\]

So,
\[\chi(\stc{T}{m+2})=f_0+f_2-f_1=(3m^2-11m+12)\frac{(m+2)!}{24}.\]

\vspace{1ex}

\nin
{\it Example 3.}
Switching to considering the small values of $m$. Set $m:=2$, so $T$ is just an edge.
The complex $\stc{T}{n}$ is a cubical subdivision of the
$(n-2)$-dimensional sphere. $\stc{T}{3}$ is a hexagon.  $\stc{T}{4}$
is a rhombic dodecahedron, whose $f$-vector is $(14,24,12)$.

In general the $f$-vector of $\stc Tn$, when $T$ is a single edge, is
$(f_0,\dots,f_{n-2})$, where
\[f_k=\binom{n}{k}(2^{n-k}-2),\text{ for }k=0,\dots,n-2.\]

\myfig{ex1}{Examples of $\stc{T}{n}$, when $T$ is an edge}

The cubical complex $\stc Tn$ can be obtained by starting with an
$n$-cube $K$ and then deleting two opposite vertices $a$ and $b$,
together with all smaller cubes in $K$ containing $a$ or $b$. The
author is not aware whether there exists some established terminology
for these complexes, beyond the cases $n=3$ and $n=4$.

\mynewpage

\section{The topology of the Stirling complexes}

\subsection{The formulation of the main theorem}

Somewhat surprisingly, our main theorem implies that the homotopy type
of the Stirling complexes $\stc Tn$ only depends on $n$ and on the
number of vertices in $T$, not on the actual tree structure.

%We can now state our main result, which completely describes the topology
%of the Stirling complexes $\stc{T}{n}$.

\begin{thm} \label{thm:main}
Assume $T$ is a tree with $m$ vertices, $m\geq 2$, and $n$ is an
integer, $n\geq m$.  The cubical complex $\stc{T}{n}$ is homotopy
equivalent to a wedge of $(n-m)$-dimensional spheres.

Let $f(m,n)$ denote the number of these spheres. Then $f(m,n)$ is
given by the following formula
\begin{multline} \label{eq:fmn}
f(m,n)=(m-1)^n-\binom{m}{1}(m-2)^n+\binom{m}{2}(m-3)^n+\dots \\
+(-1)^{m-1}\binom{m}{m-3}2^n+(-1)^m\binom{m}{m-2}.
\end{multline}
\end{thm}

In particular, we have  $f(2,n)=1$, confirming our observation that $\stc Tn$ is a sphere in this case.
Further values of $f(-,-)$ are
\begin{align*}
f(3,n) &=2^n-3, \\
f(4,n) &=3^n-4\cdot 2^n+6.
\end{align*}

\Cref{table:fmn} gives the values of $f(m,n)$ for small $m$ and $n$.
\begin{table}[hbt]
\[\begin{tabular}{ c|cccccc }
$n\setminus m$ & 2 & 3 & 4 & 5 & 6 & \dots \\  [2pt]
 \hline \\ [-5pt]
 2 & 1 &  &  &  & & \\ [2pt]
 3 & 1 & 5 &  &  & &\\   [2pt]
 4 & 1 & 13 & 23 & & &  \\  [2pt]
 5 & 1 & 29 & 121 & 119 & & \\ [2pt]
 6 & 1 & 61 & 479 & 1083 & 719 & \\  %[2pt]
$\cdots$ & $\cdots$ & $\cdots$ & $\cdots$ & $\cdots$ & $\cdots$ & $\ddots$
\end{tabular} \]
\caption{The values of $f(m,n)$.}
\label{table:fmn}
\end{table}

%\begin{rem}
It is interesting to think about the implications of \Cref{thm:main}
for the original problem of resource distribution. Clearly, the fact
that $\stc Tn$ is connected, when $n>m$, means that, starting from any
distribution one can get to any other one by moving the resources.
When $n>m+1$, the space $\stc Tn$ is simply connected. This means that
when two distributions are fixed, any two redistribution schemes from
the first distribution to the second one, are homotopic, i.e., there
is a simultaneous redistribution scheme, connecting the two. Even
higher connectivity of $\stc Tn$ means the presence of these
higher-dimensional redistribution schemes. Finally, the fact that the
homotopy of $\stc Tn$ is not trivial in the top dimension means that
in this dimension there is a number of fundamentally different
higher-dimensional redistribution schemes. The number $f(m,n)$ tells
us, in a certain sense, just how many of these different schemes there
are.
%\end{rem}

Let us make a few comments on the numerical side. First, by
Euler-Poincar\'e formula, \Cref{eq:chistc} could be used
instead of \Cref{eq:fmn}, although the latter is clearly simpler.

Second, let $\srf{m}{n}$ denote the number of surjective functions
from $[n]$ to $[m]$. We have $\srf{m}{n}=m!\sn mn$. We can
then rewrite \Cref{eq:fmn} as follows.

\begin{prop}
For all $n\geq m\geq 2$, we have
\begin{equation} \label{eq:chistar}
f(m,n)=\srf{m-1}{n}-\srf{m-2}{n}+\dots+(-1)^m\srf{1}{n}.
\end{equation}
\end{prop}
\pr
As a simple corollary of the principle of inclusion-exclusion we have the following
well-known formula
\begin{equation} \label{eq:sf}
\srf ab=a^b-\binom{a}{a-1}(a-1)^b+\dots+(-1)^{a-1}\binom{a}{1}.
\end{equation}
Substituting the right hand side of \Cref{eq:sf} into \Cref{eq:chistar},
and using the Pascal triangle addition rule for the binomial coefficients,
shows that it is equivalent to \Cref{eq:fmn}.
\qed

Finally, for future reference, we record the following fact.

\begin{prop} \label{prop:mfact}
We have \[\sum_{\alpha=1}^{m-1}(-1)^{m+\alpha+1}\binom{m}{\alpha+1}\alpha^m=m!-1.\]
\end{prop}
\pr
This follows from the following well-known polynomial identity
\[m!=\sum_{k=0}^m\binom{m}{k}(-1)^k(x-k)^m,\]
where $x$ is a variable, by substituting $x:=m+1$.
\qed

\nin
\Cref{prop:mfact} shows that \Cref{eq:fmn} holds for $m=n$.

\mynewpage

\subsection{Relaxing the occupancy requirement}

Our proof of \Cref{thm:main} proceeds by induction. As it often happens
in such situations, it is opportune to deal with a more general class of complexes.
In our case, we relax the requirement that each node must have at least one allocated
resource.

\begin{df}
\label{df:str}
For any cell $c$ of $T^n$, we define
\[\supp c=\myset{v\in V(T)}{\exists k\in[n],\textrm{ such that }c_k=v}\subseteq V(T).\]

Let $S$ be an arbitrary subset of the vertex set of $T$. We define
$\stm{T}{S}{n}$ to be the subcomplex of $T^n$, consisting of all cells
$c$ whose support contains $S$.
\end{df}

Note, that whenever $b,c\in T^n$ are cubes, such that $b\subseteq c$,
we have $\supp c\subseteq\supp b$. In other words, the support of a cell
$c$ either stays the same or increases when taking the boundary. This implies
that the cubical complex $\stm TSn$ is well-defined.

Extreme values for $S$ give us two special cases:
\begin{itemize}
\item for $S=V(T)$, we have $\stm TSn=\stc Tn$;
\item for $S=\emptyset$, we have $\stm TSn=T^n$, which is contractible
as a topological space.
\end{itemize}

Rather than attacking \Cref{thm:main} directly, we shall prove the following,
more general result.

\begin{thm} \label{thm:main2}
The complex $\stm TSn$ is homotopy equivalent to a wedge of $(n-|S|)$-dimensional
spheres. The number of spheres is $f(|S|,n)$.
\end{thm}

Clearly, \Cref{thm:main} is a special case of \Cref{thm:main2}, where
$S=V(T)$.

%\mynewpage

\section{Homotopy colimits}

\subsection{The diagrams of topological spaces}

Our strategy to prove \Cref{thm:main2} is to decompose the spaces
$\stm TSn$ into simpler pieces and then to manipulate this
decomposition, while preserving the homotopy type of the total space.
Although there are different ways to formulate our argument, we find it
handy to phrase it using the language of homotopy colimits.  Let us
introduce the corresponding terminology, see also \cite{BoK,Hat,V}.

We assume that the reader is familiar with basic category theory, \cite{ML,Mi}.
Recall, that given a poset $P$, we can always view $P$ as a category
so that
\begin{itemize}
\item the objects of that category are precisely the elements of $P$;
\item for any two elements $p,q\in P$, such that $p\geq q$, there
  exists a \emph{unique} morphism from $p$ to $q$.
\end{itemize}
The composition rule in this category is clearly uniquely defined, since there
is at most one morphism between any two objects.

Recall that $\ctop$ denotes the category of topological spaces and
continuous maps.

\begin{df}
Assume, we are given a poset $P$, and we view it as a category.
A functor from $P$ to $\ctop$ is called a \demph{diagram of topological spaces over $P$}.
\end{df}

Specifically, a diagram $\cld$ is a collection of topological spaces $\cld(p)$, where $p\in P$,
together with continuous maps $\cld_{p,q}:\cld(p)\ra\cld(q)$, where $p>q$.
These maps are subject to the condition $\cld_{q,r}\circ\cld_{p,q}=\cld_{p,r}$, whenever $p>q>r$.

%\mynewpage

\subsection{Homotopy colimits of diagrams over $P^T$}

Let $T$ be a an arbitrary tree with $m$ vertices, where $m\geq 2$. We
assume that the vertices are indexed by the set
$[m]=\{1,\dots,m\}$. A~poset $P^T$ is defined as follows:
\begin{itemize}
\item the elements of $P^T$ are indexed by the vertices and the edges
  of~$T$;
\item the partial order on $P^T$ given by saying that each edge is larger
  than its adjacent vertices.
\end{itemize}
This poset has $2m-1$ elements. The elements indexed by the vertices
are minimal, while the elements indexed by the edges are maximal, and
each one is larger than exactly $2$ minimal elements.

A diagram  $\cld$ of topological spaces over $P^T$ is then given by the following data,
subject to no further conditions:
\begin{itemize}
\item spaces $\cld(v)$ for all vertices of $T$;
\item spaces $\cld(e)$ for all edges of $T$;
\item continuous maps $\cld_{e,v}:\cld(e)\ra\cld(v)$, whenever $v$ is a vertex adjacent
to the edge $e$.
\end{itemize}

\begin{df} \label{df:hocolim}
Assume $\cld$ is a diagram of topological spaces over a poset $P^T$.
We define the {\demph homotopy colimit} of $\cld$, denoted $\hocolim\cld$,
as the quotient space
\[\hocolim\cld=\left(\coprod_{v\in V(T)} \cld(v)
\coprod_{e\in E(T)}(\cld(e)\times[0,1])\right)/\sim,\]
where the equivalence relation $\sim$ is generated by
$(x,0)\sim \cld_{e,v}(x)$,  and $(x,1)\sim \cld_{e,w}(x)$, whenever $x\in\cld(e)$,
$e=(v,w)$, $v<w$.
\end{df}

Let us mention that the notion of homotopy colimit can be defined more
generally, including homotopy colimits of diagrams of topological spaces over
arbitrary posets. Here, we restrict ourselves to \Cref{df:hocolim}, which will
be sufficient for our purposes.

%\mynewpage

\subsection{Homotopy independence of the homotopy colimits of diagrams of CW complexes
over $P^T$}

From now on, we assume that the spaces $\cld(p)$ are CW complexes, for
all $p\in P^T$, and the maps $\cld_{e,v}$ are cellular.
The next proposition says that changing these maps up to homotopy does not
change the homotopy type of the homotopy colimit.

\begin{prop} \label{prop:hhl}
Assume $\cld$ and $\cle$ are diagrams of CW complexes over $P^T$, such that
\begin{enumerate}
\item[(1)] $\cld(p)=\cle(p)$, for all $p\in P^T$;
\item[(2)] the maps $\cld_{e,v}$ and  $\cle_{e,v}$ are homotopic, whenever $e$ is an edge
of $T$, and $v$ is a vertex adjacent to $e$.
\end{enumerate}
Then $\hocolim \cld$ and $\hocolim\cle$ are homotopy equivalent.
\end{prop}
\pr Since $T$ is finite, it is enough to consider the case where
$\cld_{e,v}$ and $\cle_{e,v}$ coincide, for all, but one single
instance of an edge $e$ and a vertex $v$.

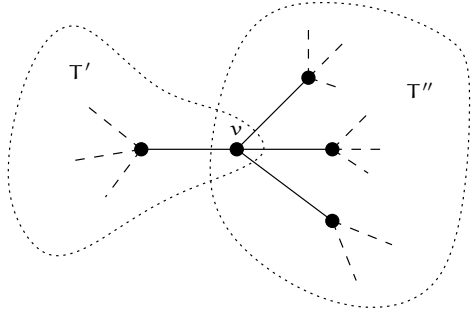
\myfig{treed}{Decomposition of the tree $T$}

Decompose the tree $T$ into a union of trees $T'$ and $T''$, such that
the intersection of $T'$ and $T''$ is vertex $v$, $v$ is a leaf of
$T'$, and $T'$ contains the edge $e$, see \Cref{fig:treed}.  Let
$\cld'$ be the diagram of CW complexes on $P^{T'}$, which is a
restriction of $\cld$ with a slight change at $v$. Specifically, it is
defined as follows:
\begin{itemize}
\item for any vertex $w\in V(T')$, we have
$\cld'(w):=\begin{cases}
\cld(w),&\text{ if } w\neq v; \\
\cld(e),&\text{ otherwise.}
\end{cases}$
\item $\cld'(r)=\cld(r)$, for all $r\in E(T')$;
\item for any edge $r\in E(T')$ and an adjacent vertex $w$, we have
\[\cld'_{r,w}=\begin{cases}
\cld_{r,w},&\text{ if } (r,w)\neq (e,v);\\
\id_{\cld(e)},&\text{ otherwise.}
\end{cases}.
\]
\end{itemize}
Let $\cld''$ be the restriction of $\cld$ to $P^{T''}$.

Set $X:=\hocolim\cld'$, $Y:=\hocolim\cld''$, $A:=\cld'(v)=\cld(e)$.
Note that $X$ and $Y$ are CW complexes, and $A$ is a CW subcomplex of
$X$.  Set $f:=\cld_{e,v}$ and $g:=\cle_{e,v}$.  Clearly,
$\hocolim\cld$ is obtained from $Y$ by attaching $X$ over $f$, whereas
$\hocolim\cle$ is obtained from $Y$ by attaching $X$ over $g$. We
assumed that $f$ is homotopic to $G$. It is then a general fact, see
e.g. \cite{Hat}, that the homotopy type of the adjunction space does
not change, when the attachment map is replaced by a homotopic
one. This implies that $\hocolim \cld$ and $\hocolim\cle$ are homotopy
equivalent.  \qed

%\mynewpage
\subsection{Special homotopy colimits}

As above, let $T$ be an arbitrary tree with at least $2$ vertices.
Let us fix a nonempty subset $S\subseteq V(T)$.  Assume we have a
diagram of CW complexes over $P^T$ satisfying the following
conditions:
\begin{itemize}
\item $\cld(v)$ are single points, for all $v\in S$;
\item $\cld(e)=X$, for all $e\in E(T)$, and $\cld(v)=X$,  for any $v\notin S$,
where $X$ is some fixed CW complex;
\item the maps $\cld_{e,v}$ are identity maps, for all $v\notin S$.
\end{itemize}

\begin{prop} \label{prop:gluing1}
Under the conditions above, the homotopy colimit of $\cld$
is homotopy equivalent to the wedge of $|S|-1$ copies of $\susp X$.
\end{prop}
\pr The proof is by induction on the number of vertices of $T$.

The induction base is when $m=2$. We have two cases.

\nin {\bf Case 1.}
If $S=1$, then $\hocolim\cld$ is a cone over $X$, hence contractible.

\nin {\bf Case 2.} If $S=2$, then  $\hocolim\cld$ is obtained by taking a cylinder over $X$ and
shrinking each of the end copies of $X$ to a point. This is precisely the suspension space $\susp X$.

From now on we can assume $m\geq 3$. We break our argument in the following cases.

\nin {\bf Case 1.}
Assume there exists an internal vertex $v\in T$, such that $v\in S$.
Let $e_1,\dots,e_k$ be the edges adjacent to $v$, $k\geq 2$.

Cutting $T$ at $v$ will decompose $T$ into the trees $T_1,\dots,T_k$,
$v$ is a leaf in each of them, and $e_i$ is adjacent to $v$ in $T_i$,
for all $i=1,\dots,k$; see \Cref{fig:treed2}. Let $\cld_i$ be the restriction of $\cld$ to
$T_i$, for $i=1,\dots,k$.

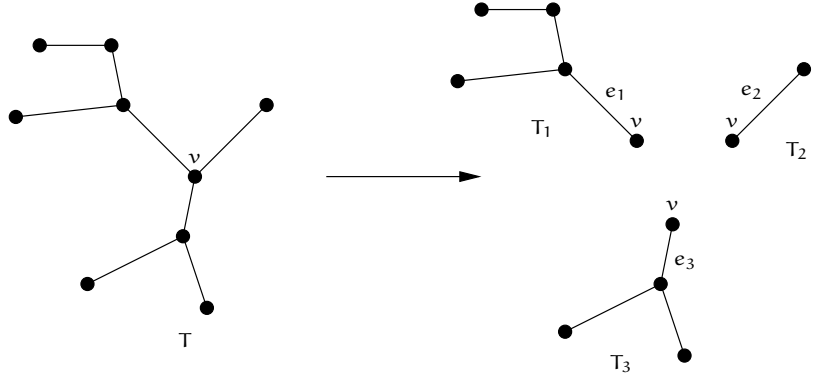
\myfig{treed2}{Cutting the tree $T$ at $v$}

Each homotopy colimit $\hocolim\cld_i$ has a marked point $x_i$
corresponding to the copy of $\cld(v)$. The homotopy colimit
$\hocolim\cld$ is obtained by gluing the homotopy colimits
$\hocolim\cld_i$ together along these points, for
$i=1,\dots,k$. Accordingly, we see that
\begin{equation} \label{eq:gl1}
\hocolim\cld\cong\vee_{i=1}^k\hocolim\cld_i.
\end{equation}

Set $S_i:=S\cap V(T_i)$. The vertex $v$ is in $S$, so $v\in S_i$, for
all $i$. This means that $S\sm v=\coprod_{i=1}^k(S_i\sm v)$, and hence
$|S|-1=\sum_{i=1}^k(|S_i|-1)$.  Since each $T_i$ has fewer vertices
than $T$, we know by the induction assumption that $\hocolim\cld_i$ is
homotopy equivalent to $|S_i|-1$ copies of $\susp X$. Accordingly,
\eqref{eq:gl1} implies that $\hocolim\cld$ is homotopy equivalent to a
wedge of $|S|-1$ copies of $\susp X$.

\nin {\bf Case 2.} All the vertices in $S$ are leaves of $T$, and
there exists a further leaf $w\notin S$.

\nin Assume $w$ is connected to the vertex $u$. Since $m\geq 3$ and
all the vertices in $S$ are leaves, we must have $u\notin S$. Let $T'$
be the tree obtained from $T$ by deleting $w$ and the adjacent edge.
Let $\cld'$ be the restriction of $\cld$ to $T'$. By induction
assumption $\hocolim\cld'$ is homotopy equivalent to a wedge of
$|S|-1$ copies of $\susp X$.  The space $\hocolim\cld$ is obtained
from $\hocolim\cld'$ by attaching a cylinder with base $X$ at one of
its ends.  Clearly $\hocolim\cld'$ is a strong deformation of
$\hocolim\cld$, so the latter is also homotopy equivalent to a wedge
of $|S|-1$ copies of $\susp X$.

\nin {\bf Case 3.} The set $S$ is precisely the set of all leaves of $T$.

\nin Since $m\geq 3$, we have at least $3$ leaves. Fix $v\in S$. Say
$v$ is connected to $w$ by an edge.  We have $w\notin S$. Let $T'$ be
the tree obtained from $T$ by deleting $v$, and let $\cld'$ be the
restriction of $\cld$ to $T'$. The topological space $\hocolim\cld$ is
obtained from $\hocolim\cld'$ by attaching a cone over
$X=\cld(w)$. Let $u\in S$ be any other leaf of $T$, $u\neq v$.  There
is a unique path inside of $T'$ connecting $w$ with $u$. The homotopy
colimit of the restriction of $\cld$ to that path is a cone with apex
at $\cld(u)$ and base at $\cld(w)$. This cone lies inside $\hocolim
T'$, therefore the inclusion map $\cld(w)\hookrightarrow\hocolim\cld'$
is trivial.  It follows that, up to homotopy equivalence, attaching a
cone over $\cld(w)$ to $\hocolim\cld'$ is the same as wedging
$\hocolim\cld'$ with $\susp X$. The result now follows by induction.
\qed

Let us now consider a little more general diagrams. These satisfy the
same conditions outside of $S$, however, for any $v\in S$, the spaces
$\cld(v)$ are now arbitrary connected CW complexes, and each
$\cld_{e,v}$ maps everything to some point in $\cld(v)$, whenever $e$
is an adjacent edge. In this case, \Cref{prop:gluing1} can be
generalized as follows.

\begin{prop} \label{prop:gluing2}
Under the conditions above the homotopy colimit of $\cld$
is homotopy equivalent to the wedge
\[\vee_{v\in S}\cld(v)\vee_\Omega\susp X,\]
where $|\Omega|=|S|-1$.
\end{prop}
\pr For each $v\in V(T)$ we select a base point $x_v\in\cld(v)$. Since
$\cld(v)$ is connected, any continuous map $\varphi:Y\ra\cld(v)$,
mapping everything to a point, is homotopic to a map
$\psi:Y\ra\cld(v)$ mapping everything to the base point $x_v$. By
\Cref{prop:hhl} we can therefore assume that each map $\cld_{e,v}$
maps everything to $x_v$, whenever $v\in S$ and $e$ is an adjacent
edge, without changing the homotopy type of $\hocolim\cld$.

Let $\cld'$ be the diagram which we obtain from $\cld$ by replacing
each $\cld(v)$ by a point, for $v\in S$.  Clearly, the homotopy
colimit $\hocolim\cld$ is obtained from $\hocolim\cld'$ be wedging it
with all the $\cld(v)$, for $v\in S$. The result follows from
\Cref{prop:gluing1}.  \qed

%\mynewpage

\section{Structural decomposition of Stirling complexes and consequences}

\subsection{Representing Stirling complexes as homotopy colimits.} $\,$

\nin
Let us fix $n$ and $S$.
We define a diagram $\cld$ of topological spaces over $P^T$ as follows:
\begin{itemize}
\item for each vertex $v\in V(T)$, we set $\cld(v)$ to be the subcomplex of $\stm TSn$,
consisting of all cells $c$, such that $c_n=v$;
\item for each edge $e\in E(T)$, we set $\cld(e)$ to be the Stirling
complex $\stm{T}{S}{n-1}$;
\item finally, for each edge $e\in E(T)$, $e=(v,w)$, we define the map
$\cld_{e,v}:\cld(e)\ra\cld(v)$ by setting $c_n:=v$.
\end{itemize}

\myfig{hocolim}{$\stm{T}{S}{n}$ as $\hocolim\cld$}

\begin{prop} \label{prop:repr}
The homotopy colimit $\hocolim\cld$ is homeomorphic to
the cubical complex $\stm{T}{S}{n}$.
\end{prop}
\pr Whenever $e=(v,w)$ is an edge of $T$, let $B_e$ denote
the subcomplex of $\stm TSn$ consisting of all cells $c\in\stm TSn$ such that
one of the following holds
\begin{enumerate}
\item[(1)] $c_n=e$;
\item[(2)] $c_n=v$, and there exists $1\leq k\leq n-1$, such that
$c_k=v$;
\item[(3)] $c_n=w$, and there exists $1\leq k\leq n-1$, such that
$c_k=w$.
\end{enumerate}
It is easy to see that this set of cells is closed under taking the
boundary, hence the subcomplex $B_e$ is well-defined.  Furthermore,
the complex $\stm TSn$ is the union of the subcomplexes $\cld(v)$, for
$v\in V(T)$, and $B_e$, for $e\in E(T)$. To see this just take any
cube $(c_1,\dots,c_n)$ and sort it according to the value of $c_n$.

Recording the value of $c_n$ separately, we can see that,
as a cubical complex, each $B_e$ is isomorphic to the direct product of
 $\stm TS{n-1}$ with the closed interval $[0,1]$. This can be seen as a cylinder
with base $\stm TS{n-1}$. The entire complex
$\stm TSn$ is obtained by taking the disjoint union of  $\cld(v)$, for $v\in V(T)$,
and connecting them by these cylinders. For each cylinder $B_e$, $e=(v,w)$, its bases are identified
with corresponding subcomplexes of $\cld(v)$ and $\cld(w)$ by assigning
$c_n:=v$ or $c_n:=w$. These are precisely the maps $\cld_{e,v}$ and $\cld_{e,w}$.
Comparing this gluing procedure with the definition of $\hocolim\cld$ we see that
we obtain a homeomorphic space.
\qed

\subsection{The proof of the main theorem}

We are now ready to show our main result.

\vspace{1ex}

\prn{\Cref{thm:main2}} First, when $|S|=n$, the complex $\stm TSn$ is a disjoint union
of $n!$ points. This can be viewed as a wedge of $n!-1$ copies of a $0$-dimensional sphere,
so the result follows from \Cref{prop:mfact}.

Assume from now on that $n\geq |S|+1$.  By \Cref{prop:repr} we can
replace $\stm TSn$ by $\hocolim\cld$.  Consider now a map
$\cld_{e,v}:\cld(e)\ra\cld(v)$. By induction, we know that
$\cld(e)=\stm{T}{S}{n-1}$ is homotopy equivalent to a wedge of spheres
of dimension $n-1-|S|$. We make 2 observations.
\begin{enumerate}
\item[(1)] If $v\notin S$, the cubical complex $\cld(v)$ is isomorphic to $\stm TS{n-1}$,
and the map $\cld_{e,v}$ is the identity map.
\item[(2)] If $v\in S$, the cubical complex $\cld(v)$ is isomorphic to $\stm T{S\sm v}{n-1}$.
This is because we know that $c_n=v$, so there is no need to request that $v$ is occupied by
some other resource. By induction assumption, the space  $\stm T{S\sm v}{n-1}$ is homotopy
equivalent to a wedge of spheres of dimension $n-1-(|S|-1)=n-|S|$. In particular, it is
$(n-|S|-1)$-connected. Therefore, the map $\cld_{e,v}$ is homotopic to a trivial map, which
takes everything to a point.
\end{enumerate}

We now apply \Cref{prop:hhl} to shift our consideration to the diagram
$\cld'$, which is obtained from $\cld$ by replacing the maps
$\cld_{e,v}$ by trivial ones, whenever $v\in S$. This diagram has the
same homotopy type as $\cld$. On the other hand, it now satisfies the
conditions of \Cref{prop:gluing2}, where the connectivity of the
spaces $\cld(v)$ is a consequence of the fact that $n\geq|S|+1$. It
follows from that proposition that
\[\hocolim\cld\simeq \vee_{v\in S}\cld(v)\vee_\Omega\susp\stm TS{n-1},\]
where $|\Omega|=|S|-1$.

Counting spheres on both sides, we obtain the recursive formula
\[f(|S|,n)=(|S|-1)f(|S|,n-1)+|S|f(|S|-1,n-1).\]
The validity of the formula \Cref{eq:fmn} now follows from \Cref{prop:fmn}.
\qed

\begin{rem}
After this paper was submitted for publication, a shorter proof of
\cref{thm:main2} was found by one of the referees. It is included in
the appendix.
\end{rem}

%\mynewpage

\begin{prop} \label{prop:fmn}
Let $\Gamma=\myset{(m,n)\in\zz\times\zz}{n\geq m\geq 2}$.
Assume we have a function $f:\Gamma\ra\zz$, which satisfies the following:
\begin{enumerate}
\item[(1)] for all $n>m\geq 3$ we have recursive formula
\begin{equation} \label{eq:rec}
f(m,n)=(m-1)f(m,n-1)+m f(m-1,n-1);
\end{equation}
\item[(2)]
we have the boundary conditions $f(2,n)=1$, $f(m,m)=m!-1$.
\end{enumerate}

Then for all $(m,n)\in\Gamma$, the value $f(m,n)$ is given by \Cref{eq:fmn}, which we rewrite as
\begin{equation}
\label{eq:fmn2}
f(m,n)=\sum_{\alpha=1}^{m-1}(-1)^{m+\alpha+1}\binom{m}{\alpha+1}\alpha^n.
\end{equation}
\end{prop}
\pr Clearly, the recursive rule \Cref{eq:rec} together with the boundary conditions
defines the values of the function of the function $f(-,-)$ uniquely.
Therefore, to show that $f$ is given by the formula   \Cref{eq:fmn2} we just need to know that
this formula satisfies our boundary conditions and the recursion.

Substituting $m=2$ into \Cref{eq:fmn2} immediately yields $1$ on the
right hand side, as there is only one summand, with $\alpha=1$.
The case $m=n$ follows from \Cref{prop:mfact}.

To show that \Cref{eq:fmn} satisfies the recursion \Cref{eq:rec} we need to check that
\begin{multline}\label{eq:recal}
\sum_{\alpha=1}^{m-1}(-1)^{m+\alpha+1}\binom{m}{\alpha+1}\alpha^n=\\
(m-1)\sum_{\alpha=1}^{m-1}(-1)^{m+\alpha+1}\binom{m}{\alpha+1}\alpha^{n-1}+
m\sum_{\alpha=1}^{m-2}(-1)^{m+\alpha}\binom{m-1}{\alpha+1}\alpha^{n-1}.
\end{multline}
We do that simply by comparing the coefficients of $\alpha^{n-1}$ on each side of \Cref{eq:recal}.
For $\alpha=m-1$, the coefficient on each side is $m-1$.
For $\alpha=1,\dots,m-2$, we need to show that
\[(-1)^{m+\alpha+1}\binom{m}{\alpha+1}\alpha=(m-1)(-1)^{m+\alpha+1}\binom{m}{\alpha+1}
+m(-1)^{m+\alpha}\binom{m-1}{\alpha+1}.\]
This follows from the formula
\[(m-\alpha-1)\binom{m}{\alpha+1}=m\binom{m-1}{\alpha+1}. \qed\]

We finish with an open question.
\begin{openq}
Let $T$ be a tree with a single internal vertex of valency $r$, where
$r\geq 2$, and let $n$ be any integer, $n\geq r+1$. The symmetric
group $\cls_r$ acts on $T$ by permuting its $r$ leaves. This induces an
$\cls_r$-action on the Stirling complex $\stc T n$, and hence also an
$\cls_r$-action on $H_{n-r-1}(\stc Tn;\rr)$. It would be interesting
to decompose this representation of $\cls_r$ into irreducible ones.
\end{openq}

\appendix
\setcounter{equation}{0}
\renewcommand{\theequation}{A\thesection.\arabic{equation}}
\setcounter{theorem}{0}
\renewcommand{\thetheorem}{A\thesection.\arabic{theorem}}

\section*{Appendix: Stirling Complexes via the Wedge Lemma}

\smallskip
\begin{center}
by \textsc{Roy Meshulam \footnote[1]{Department of Mathematics,
Technion, Haifa 32000, Israel. e-mail:
meshulam@technion.ac.il~. Supported by ISF grant 686/20.}}
\end{center}
\medskip

In this appendix we prove a generalization of \Cref{thm:main2}.
Let $X=X_0$ be a finite simplicial complex, and let
$S=\{X_i\}_{i \in [m]}$
be a family of subcomplexes of $X$.
For $n \geq m$
let
\[
A_{m,n}=\left\{ (i_1,\ldots,i_n) \in (\{0\}\cup [m])^n : \{i_1,\ldots,i_n\} \supset [m] \right\}.
\]
Slightly extending the setup considered in \ref{df:str}, we define the \emph{Stirling Complex} associated with the triple $(X,S,n)$ by
\[
\str(X,S,n)= \bigcup_{(i_1,\ldots,i_n) \in A_{m,n}} X_{i_1} \times \cdots \times X_{i_{n}}.
\]
%For $1 \leq m \leq n$ define $f(m,n)$ inductively by $f(1,n)=0$, $f(m,m)=m!-1$ and
%\begin{equation}
%\label{e:recf}
%f(m,n)=mf(m-1,n-1)+(m-1)f(m,n-1)
%\end{equation}
%for $n>m \geq 2$. 
Let $\sph^k$ denote the $k$-sphere. 
\Cref{thm:main2} asserts that 
if $T$ be a finite tree and $S$ is a set of $m \geq 2$ distinct vertices of $T$, then
\[
\str(T,S,n)\simeq \bigvee_{i=1}^{f(m,n)} \sph^{n-m}.
\]
\noindent
Here we give a simple proof of a generalization of \Cref{thm:main2} that perhaps clarifies why
the homotopy type of $\str(T,S,n)$ does not depend on the structure of $T$.

\begin{theorem}
\label{t:app}
Let $X$ be a finite contractible complex and let $S=\{X_i\}_{i=1}^m$ be a family of $m \geq 2$ pairwise disjoint contractible subcomplexes of $X$. Then
\[
\str(X,S,n)\simeq \bigvee_{i=1}^{f(m,n)} \sph^{n-m}.
\]
\end{theorem}
\noindent
The main tool in the proof of Theorem \ref{t:app} is the Wedge Lemma of
Ziegler and \v{Z}ivaljevi\'{c} (Lemma 1.8 in \cite{ZZ}). The version below appears in
\cite{HRW}. For a poset $(P,\prec)$ and $p \in P$ let $P_{\prec p}=\{q \in P: q \prec p\}$. Let $\Delta(P)$ denote the order complex of $P$. Let $Y$ be a regular CW-complex and let $\{Y_i\}_{i=1}^m$ be subcomplexes of $Y$ such that $\bigcup_{i=1}^m Y_i=Y$. Let $(P,\prec)$ be the poset whose elements index all distinct partial intersections $\bigcap_{j \in J} Y_j$, where $\emptyset \neq J \subset [m]$. Let $U_p$ denote the partial intersection indexed by $p \in P$,
and let $\prec$ denote reverse inclusion, i.e. $p \prec q$ if $U_q \subsetneqq U_p$.

\vspace{2mm}
\noindent
{\bf Wedge Lemma \cite{ZZ,HRW}.}
suppose that for any $p \in P$ there exists a $c_p \in U_p$ such that the inclusion
$\bigcup_{q \succ p} U_q \hookrightarrow U_p$ is homotopic to the constant map to $c_p$. Then
\begin{equation}
\label{e:wlemma}
Y \simeq \bigvee_{p \in P} \Delta(P_{\prec p})*U_p.
\end{equation}
\noindent
{\bf Proof of Theorem \ref{t:app}.} If $m=n$ then $\str(X,S,n)$ is a union of $m!$ disjoint contractible sets and hence homotopy equivalent to $\bigvee_{i=1}^{m!-1} \sph^0$. Suppose $n>m \geq 2$.
In view of the recursion (\ref{eq:rec}), it suffices as in the proof of \Cref{thm:main2}, to establish the following homotopy decomposition:
\begin{equation}
\label{e:homdec}
\str(X,S,n) \simeq \bigvee_{i=1}^m \str(X,S \setminus\{ X_i \},n-1) \vee \bigvee_{i=1}^{m-1} \sph^0  * \str(X,S,n-1) .
\end{equation}
We proceed with the proof of (\ref{e:homdec}).
For $1 \leq i \leq m$ let
$$
Y_i=\big(X_i \times \str(X,S \setminus \{X_i\},n-1)\big) \cup
\big(X \times \str(X,S,n-1) \big).
$$
Then $\bigcup_{i=1}^m Y_i=\str(X,S,n)$. Next note that
\begin{equation}
\label{e:xsubs}
\str(X,S,n-1) \subset \str(X,S \setminus \{X_i\},n-1).
\end{equation}
As $X_i \subset X$ are both contractible, it follows that $X_i$ is a deformation retract of $X$.
Together with (\ref{e:xsubs}) it follows that
$X_i \times \str(X,S \setminus \{X_i\},n-1)$ is a deformation retract of $Y_i$. Therefore
\begin{equation}
\label{e:yv}
Y_i \simeq \str(X,S \setminus \{X_i\},n-1).
\end{equation}
Let $Z=X \times \str(X,S,n-1)$.
Then for any $1 \leq i \neq j \leq m$
\begin{equation}
\label{e:capyv}
Y_{i} \cap Y_{j}= \bigcap_{k=1}^m  Y_k = Z \simeq \str(X,S,n-1).
\end{equation}
%and
%\begin{equation}
%\label{e:zz}
%Z \simeq \str(X,S,n-1).
%\end{equation}
\noindent
Eq. (\ref{e:capyv}) implies that the intersection poset  $(P,\prec)$ of the cover $\{Y_i\}_{i=1}^m$ is $P=[m] \cup \{\widehat{1}\}$, where $i \in [m]$ represents $Y_i$, $\widehat{1}$ represents $Z$,
$[m]$ is an antichain and
$i \prec \widehat{1}$ for all $i \in [m]$.
Note that
$\Delta(P_{\prec i})=\emptyset$ for all $i \in [m]$ and $\Delta(P_{\prec \widehat{1}})$ is the discrete space $[m]$.
By induction, $Y_i$ is homotopy equivalent to a wedge of $(n-m)$-spheres and $Z$ is homotopy equivalent to a wedge of $(n-m-1)$-spheres. Hence
the inclusion $Z \hookrightarrow Y_i$ is null homotopic. Using the Wedge Lemma together with
 (\ref{e:yv}) and (\ref{e:capyv}), it follows that
\begin{equation*}
\label{e:wedgel}
\begin{split}
\str(X,S,n) &\simeq \left( \bigvee_{i \in [m]} \Delta(P_{\prec i}) * Y_i\right) \vee
\left( \Delta(P_{\prec \widehat{1}}) * Z\right)
=\left(\bigvee_{i \in [m]}Y_i\right) \vee ([m]*Z) \\
&\simeq \bigvee_{i \in [m]}
\str(X,S \setminus \{X_i\},n-1)
 \vee \bigvee_{i=1}^{m-1} \sph^0  * \str(X,S,n-1).
\end{split}
\end{equation*}
This completes the proof of (\ref{e:homdec}) and hence of Theorem \ref{t:app}.
{\qed}


\begin{thebibliography}{AaA00}

\bibitem[BaK06]{BaK06} E.\ Babson, D.N.\ Kozlov, {\em Complexes of graph
homomorphisms},  Israel J.\ Math.\ {\bf 152} (2006), 285--312.

\bibitem[BaK07]{BaK07} E.\ Babson, D.N.\ Kozlov, {\em Proof of the
    Lov\'asz Conjecture}, Annals of Math.\ (2) {\bf 165} (2007), no.\ 3,
965--1007.

\bibitem[BoK]{BoK} M.\ Bousfield, D.M.\ Kan, {\it Homotopy Limits,
Completions and Localizations}, Springer Lect.\ Notes Math.\ {\bf
304}, Berlin-Heidelberg-New York 1972.

\bibitem[Ca]{Ca} G.\ Carlsson, {\em Topology and data}, Bull.\ Amer.\ Math.\
Soc.\ (N.S.) {\bf 46} (2009), no.\ 2, 255--308.

\bibitem[FFG]{FFG} A.T.\ Fomenko, D.B.\ Fuks, V.L.\ Gutenmacher, {\em
Homotopic topology}, Translated from the Russian by
K. M\'alyusz. Akad\'emiai Kiad\'o (Publishing House of the Hungarian
Academy of Sciences), Budapest, 1986.

\bibitem[Fu]{Fu} W.\ Fulton, {\em Algebraic topology}, Graduate Texts in
Mathematics {\bf 153}, Springer-Verlag, New York, 1995. xviii+430 pp.

\bibitem[GKP]{GKP} R.L.\ Graham, D.E.\ Knuth, O.\ Patashnik,
{\it Concrete Mathematics: A Foundation for Computer Science,} 2nd
ed. Reading, MA: Addison-Wesley, pp. 257-267, 1994.

\bibitem[GH]{GH} M.J.\ Greenberg, J.R.\ Harper, {\em Algebraic Topology},
Mathematics Lecture Note Series {\bf 58}, Benjamin/Cummings Publishing
Co., Inc., Advanced Book Program, Reading, Mass., 1981. xi+311 pp.

\bibitem[Har]{Har} F.\ Harary, {\em Graph Theory}, Addison-Wesley
  Series in Mathematics, Reading, MA, 1969.

\bibitem[HRW]{HRW}
J.\ Herzog, V.\ Reiner and V.\ Welker, {\it The Koszul property in affine semigroup rings},
Pacific J.\ Math.\ {\bf 186} (1998), pp.\ 39-–65.

\bibitem[Hat]{Hat} A.\ Hatcher, {\it Algebraic topology}, Cambridge
University Press, Cambridge, 2002.

\bibitem[HKR]{HKR} M.\ Herlihy, D.N.\ Kozlov, S.\ Rajsbaum,
{\it Distributed computing through combinatorial topology},
Elsevier/Morgan Kaufmann, Waltham, MA, 2014. xiv+319 pp.

\bibitem[Kn]{Kn} D.E.\ Knuth, {\it The Art of Computer Programming, Vol. 1:
  Fundamental Algorithms,} 3rd ed. Reading, MA: Addison-Wesley, 1997.

\bibitem[Ko07]{Ko07} D.N.\ Kozlov,
{\em Chromatic numbers, morphism complexes, and Stiefel-Whitney
characteristic classes}, in: {\em Geometric Combinatorics} (eds.\ E.\
Miller, V.\ Reiner, B.\ Sturmfels), pp.\ 249--315, IAS/Park City
Mathematics Series {\bf 13}, American Mathematical Society,
Providence, RI; Institute for Advanced Study (IAS), Princeton, NJ;
2007.

\bibitem[Ko08]{Ko08} D.N.\ Kozlov, \textit{Combinatorial Algebraic Topology},
Algorithms and Computation in Mathematics {\bf 21}, Springer, Berlin,
2008, xx+389 pp.

\bibitem[Ko20]{Ko20} D.N.\ Kozlov, \textit{Organized collapse: an
  introduction to discrete Morse theory}, Graduate Studies in
  Mathematics {\bf 207}, American Mathematical Society, Providence,
  RI, 2020, xxiii+312 pp.

\bibitem[ML]{ML} S.\ MacLane, {\em Categories for the Working
    Mathematician}, Second edition, Graduate Texts in Mathematics,
  {\bf 5}, Springer-Verlag, New York, 1998.

\bibitem[Mi]{Mi} B.\ Mitchell, {\it Theory of categories},
 Pure and Applied Mathematics, Vol.\ XVII, Academic Press,
New York-London, 1965.

\bibitem[Mu]{Mu} J.R.\ Munkres, {\em Elements of algebraic topology},
Addison-Wesley Publishing Company, Menlo Park, CA, 1984 ix+454 pp.

\bibitem[S]{S} J.\ Stirling, {\em Methodus differentialis, sive tractatus de summation et
interpolation serierum infinitarium.} London, 1730. English translation
by Holliday, J. The Differential Method: A Treatise of the Summation
and Interpolation of Infinite Series. 1749.

\bibitem[V]{V} R.M.\ Vogt, {\it Homotopy limits and colimits},
Math.\ Z.\ {\bf 134} (1973), pp.\ 11--52.

\bibitem[ZZ]{ZZ}
G.M.\ Ziegler and R.\ \v{Z}ivaljevi\'{c}, {\it Homotopy types of subspace arrangements via diagrams
of spaces}, Math.\ Ann.\ {\bf 295} (1993) pp.\ 527--548.



\end{thebibliography}
\end{document}